\documentclass[final]{siamltex}
\usepackage{leftidx}
\usepackage{mathtools}
 \usepackage{mathrsfs}
\usepackage{amssymb}
 \usepackage{amsmath}
 \usepackage{cite}
 \usepackage{float}
 \usepackage{multirow}
 \usepackage{color}

\def\af{\alpha}

\def\bt{\beta}

\def\Lam{\Lambda}

\def\mp{\mathcal{P}}
\def\mr{\mathbb{R}}

\newtheorem{remark}[theorem]{Remark}
\def\zhu#1{\textcolor{black}{#1}}
 \usepackage{leftidx}

\title{Superconvergence points of fractional spectral interpolation \thanks{This work is supported in part by the National Natural Science Foundation of China under grants 11471031 and 91430216, and the US
National Science Foundation through grants DMS-1419040.} } 

\author{Xuan Zhao \thanks{Beijing Computational Science Research Center, { Beijing 100084}, P. R. China,
(\tt xuanzhao11@gmail.com).}
       \and Zhimin Zhang\thanks{Beijing Computational Science Research Center, { Beijing 100084, P. R. China,} and Department of Mathematics, Wayne State University, {  Detroit, MI 48202}, USA, (\tt zmzhang@csrc.ac.cn; zzhang@math.wayne.edu).}}

\begin{document}

\maketitle

\begin{abstract}
 We investigate superconvergence properties of the spectral interpolation involving fractional derivatives. Our interest in this superconvergence problem is, in fact, twofold: when interpolating function values, we identify the points at which fractional derivatives of the interpolant superconverge; when interpolating  fractional derivatives, we locate those points where function values of the interpolant superconverge. For the former case, we apply various Legendre polynomials as basis functions and obtain the superconvergence points, which naturally unify the superconvergence points for the first order derivative presented in [Z. Zhang, SIAM J. Numer. Anal., 50 (2012), 2966-2985], depending on orders of fractional derivatives. While for the  latter case, we utilize Petrov-Galerkin method based on generalized Jacobi functions (GJF) [S. Chen et al., arXiv: 1407. 8303v1] and locate the superconvergence points both for function values and fractional derivatives.  Numerical examples are provided to verify the analysis of superconvergence points for each case.
\end{abstract}



\pagestyle{myheadings}
\thispagestyle{plain}
\markboth{Xuan Zhao, Zhimin Zhang}{Superconvergence of fractional polynomial interpolation}
\section{Introduction}



Superconvergence of numerical methods usually happens when  the convergence rate  at some special points  is higher than the {  theoretical global rate} \cite{Dougalis+1980+SINUM,Larsen+1982+SINUM,Wahlbin+1995,Zhang+2005+MC}. The investigation of superconvergence provides a fundamental  insight for {  post-processing and adaptive algorithm design, which leads to higher accuracy and more efficient numerical methods}  \cite{Ainsworth+2000+NY,Cockburn+2003+MC}.

For classical problems, the superconvergence {  phenomenon of the} $h$-version finite element method is
{  well understood and addressed in the literature, see, e.g., \cite{Cockburn+2003+MC,Lin+2006+Beijing,Wahlbin+1995},
 and references therein. On the other hand, there {  has been only a limited study of the superconvergence phenomenon for} spectral methods \cite{Zhang+2005+MC,Zhang+2008+JSC,Zhang+2012+SINUM}. {  Specifically}, in \cite{Zhang+2005+MC,Zhang+2008+JSC},
 Zhang identified the derivative superconvergence points of
  {  the Legendre and Chebyshev spectral collocation methods for the} two-point boundary value problem. In addition, Zhang \cite{Zhang+2012+SINUM} investigated {  various polynomial interpolations and {  located} superconvergence points
  for the function value and the first derivative}.


{  In this paper, we conduct the first superconvergence study of} fractional operators.
Since fractional calculus generalizes the classical (integer-order) differentiation and integration to any order, {  the study of superconvergent}  fractional derivatives unifies the theoretical investigation} in \cite{Zhang+2012+SINUM} for the first-order derivative.


Over the last two decades, fractional differential equations (FDEs) have been demonstrated to be more effective, when modeling some complex systems in physics, finance and et al. \cite{Man+1968+SIAMR,Mark+2001+PRE,Metzler+2000+RP,Podlubny+1999,Sab+2002+EPJ,sun+2009+pa}, \zhu{compared to classical models}. {  These new models are usually derived by replacing integer-order derivatives with fractional derivatives in classical models}.
Meanwhile, there has been a growing need for the development of high-order numerical \zhu{algorithms} for solving {  FDEs}.
\zhu{Due to} the nonlocal {  definition} of fractional derivative, existing methods including finite {  difference and} finite element  methods\cite{Mus+2012+IMA,Shen+2008+IMA,Stynes+2014+IMA,Tad+2006+JCP,HWang+2013+JCP,Yang+2011+SISC,Zhao+2014+JSC,Zhou+2013+JSC} mostly  lead to {  low-order} schemes.



{  Spectral methods are promising candidates for solving FDEs since their global nature fits well with the nonlocal definition of fractional operators. Using integer-order orthogonal polynomials as basis functions, spectral methods \cite{Li+2012+FCAA,Li+2009+SINUM,Xu+2014+JCP,Zeng+2014+SINUM} really help with the alleviation of the memory cost for discretization of fractional derivatives. Furthermore, to deal with singularities, which usually appear in fractional problems, in the works  \cite{Chen+2014+ar,Huang+2014+ar,Zayernouri+2013+JCP,Zayernouri+2014+SISC} the authors design suitable bases.

A major difficulty in the investigation of superconvergence of spectral methods for fractional problems, compared with integer-order derivatives, is the nonlocality of the fractional operator and the complicated form of fractional derivatives.
The second challenge is the construction of a good basis for a spectral scheme. Given a suitable  basis, one can then begin the analysis of the approximation error in order to locate the superconvergence points.

In this paper, we interpolate function values using \zhu{Legendre} polynomials, and use the error equation to obtain superconvergent points for \zhu{Riemann-Liouville} fractional derivatives. It turns out that these superconvergence points are zeros of the fractional derivatives of the corresponding \zhu{Legendre} polynomials. We also consider the interpolation of fractional Riemann-Liouville derivatives, and apply a GJF-spectral-Petrov-Galerkin scheme \cite{Chen+2014+ar} to solve an equivalent fractional initial-value problem. We have found that superconvergence points for numerical function values and fractional Riemann-Liouville derivatives are zeros of corresponding basis functions and Gauss points, respectively.

The organization of this paper is as follows. In Section 2, we define notation and give properties of fractional derivatives and of GJFs.
The interpolations of the function values using various \zhu{Legendre} polynomials are presented in Section 3.  Moreover, we analyse the error to identify superconvergence points for fractional derivatives of the interpolant. In Section 4, we show the locations of the superconvergence points both for the numerical solution and the fractional derivatives using a GJF-spectral-Petrov-Galerkin scheme. The numerical tests for both cases are displayed at the end of Section 3 and Section 4 respectively.
We conclude with a discussion of our results and their applications in Section 5}.

\section{Preliminaries}
In this section, we present some notations and lemmas which will be used in the following sections.

\begin{definition}
The fractional integral  of order  $\mu\in(0,1)$ for function $f(x)$ is  defined as
\begin{align}
(\prescript{}{x_L}{\mathcal{I}}_{x}^{\mu}f)(x)=\frac1{\Gamma(\mu)}\int_{x_L}^x\frac{f(s)}{(x-s)^{1-\mu}}ds,\quad x>x_L,
\end{align}
\end{definition}

\begin{definition}\label{def_c}
Caputo fractional derivative of order $\mu\in(0,1)$ for function $f(x)$ is  defined as
\begin{align}
(\prescript{C}{x_L}{\mathcal{D}}_{x}^{\mu}f)(x)=\prescript{}{x_L}{\mathcal{I}}_{x}^{1-\mu}\left[\frac{d }{dx}f(x)\right]
=\frac1{\Gamma(1-\mu)}\int_{x_L}^x\frac{f^{\prime}(s)}{(x-s)^{\mu}}ds,\quad x>x_L.
\end{align}
\end{definition}

\begin{definition}
Left Riemann-Liouville fractional derivative of order  $\mu\in(0,1)$ for function $f(x)$ is  defined as
\begin{align}
(\prescript{RL}{x_L}{\mathcal{D}}_{x}^{\mu}f)(x)=\frac{d}{dx}\left[(\prescript{}{x_L}{\mathcal{I}}_{x}^{1-\mu}f)(x)\right]=\frac1{\Gamma(1-\mu)}\frac{d}{dx}\int_{x_L}^x\frac{f(s)}{(x-s)^{\mu}}ds,\quad x>x_L,
\end{align}
\end{definition}

\begin{definition}
Right Riemann-Liouville fractional derivative of order  $\mu\in(0,1)$ for function $f(x)$ is  defined as
\begin{align}
(\prescript{RL}{x}{\mathcal{D}}_{x_R}^{\mu}f)(x)=\frac{1}{\Gamma(1-\mu)}\left(-\frac{d}{dx}\right)\int_{x}^{x_R}\frac{f(s)}{(s-x)^{\mu}}ds,\quad x<x_L,
\end{align}
\end{definition}

For the convenience, if $\mu<0,$ we denote $(\prescript{RL}{x_L}{\mathcal{D}}_{x}^{\mu}f)(x)=(\prescript{}{x_L}{\mathcal{I}}_{x}^{-\mu}f)(x)$ throughout the paper.

\begin{lemma}(see\cite{Ishteva+2005+MSRJ})\label{fracLeibC}\quad
Let $p\in\mathbb{R}$ and $n-1<p<n\in \mathbb{R}.$ If $f(t)$ and $\varphi(t)$ along with all its derivatives are continuous in $[a, t].$ Under this condition the Leibniz rule for  Caputo  fractional differentiation  takes the form:
\begin{align*}
\prescript{C}{a}{\mathcal{D}}_{t}^{p}\left[\varphi(t) f(t)\right]=\sum_{k=0}^{\infty}{p \choose k}\varphi^{(k)}(t) \prescript{C}{a}{\mathcal{D}}_{t}^{p-k}f(t)-\sum_{k=0}^{n-1}\frac{t^{k-p}}{\Gamma(k+1-p)}[\varphi(t) f(t)]^{(k)}(-1).
\end{align*}
\end{lemma}

\begin{lemma}(see\cite{Podlubny+1999} Chap. 2)\label{fracLeib}\quad
Let $p\in\mathbb{R}$ and $n-1<p<n\in \mathbb{R}.$ If $f(t)$ and $\varphi(t)$ along with all its derivatives are continuous in $[a, t].$ Under this condition the Leibniz rule for Riemann-Liouville fractional differentiation  takes the form:
\begin{align*}
\prescript{RL}{a}{\mathcal{D}}_{t}^{p}\left[\varphi(t) f(t)\right]=\sum_{k=0}^{\infty}{p \choose k}\varphi^{(k)}(t) \prescript{RL}{a}{\mathcal{D}}_{t}^{p-k}f(t).
\end{align*}
\end{lemma}


We recall the definition of GJFs, which extend the parameters of classical Jacobi polynomials to a wider range .
\begin{definition}(see \cite{Chen+2014+ar}) [Generalized Jacobi functions] Define
 \begin{align}
\prescript{+}{}{J}_{n}^{(-\af,\bt)}{ (x)}:=(1-x)^\af P_n^{(\af,\bt)}(x),\quad \text{for}\; \af >-1, \;\bt\in \mr,\label{GJFsr}
\end{align}
\text{and}
\begin{align}
\prescript{-}{}{J}_{n}^{(\af,-\bt)}{ (x)}:=(1+x)^\bt P_n^{(\af,\bt)}(x),\quad \text{for}\; \af \in\mr,\; \bt>-1,\label{GJFsl}
\end{align}
for $x\in\Lambda=(-1,1)$ and $n\geqslant0.$
\end{definition}

 In what follows, we present the selected two special cases for the fractional derivatives of GJFs:
\begin{itemize}
  \item Let $\af >0, \;\bt\in \mr$ and $\;n\in\mathbb{N}_0$
\begin{align}
\prescript{RL}{x}{\mathcal{D}}_{1}^{\af}\left\{ \prescript{+}{}{J}_{n}^{(-\af,\bt)} (x)\right\}=\frac{\Gamma(n+\af+1)}{n!}{ P}_{n}^{(0,\af+\bt)}(x).\label{fracGJFr}
\end{align}
  \item Let $\bt >0, \;\af\in \mr$ and $\;n\in\mathbb{N}_0$
\begin{align}
\prescript{RL}{-1}{\mathcal{D}}_{x}^{\af}\left\{ \prescript{-}{}{J}_{n}^{(\af,-\bt)} (x)\right\}=\frac{\Gamma(n+\bt+1)}{n!}{ P}_{n}^{(\af+\bt,0)}(x).\label{fracGJFl}
\end{align}
\end{itemize}
\zhu{For convenience{ , we introduce} two sets of parameters for GJFs:}
\begin{align}
&\prescript{+}{}{\Upsilon}_{1}^{\af,\bt}:=\{(\af,\bt): \af>0,\;\bt>-1\},\\
&\prescript{+}{}{\Upsilon}_{2}^{\af,\bt}:=\{(\af,\bt): \af>0,\;-\af-1<\bt=-k\leq -1,\; k\in\mathbb{N}\}.
\end{align}
For $(\af,\bt)\in\prescript{+}{}{\Upsilon}_{1}^{\af,\bt}$, we define the finite-dimensional fractional-polynomial space:
\begin{align}
\prescript{+}{}{\mathcal{F}}_{N}^{(-\af,\bt)}(\Lambda):=\text{span}\{\prescript{+}{}{J}_{n}^{(-\af,\bt)},0
\leq n\leq N\}.
\end{align}

\section{Legendre interpolant} In this section, we interpolate a smooth function $u$ at a set of $N +1$ special points{ , which are zeros of some orthogonal polynomials on $[-1,1]$. In particular, we want to find} $u_N$ such that
\begin{align}
u_N(x_k)=u(x_k),\quad -1\leq x_0<x_1<\cdots<x_N\leq1,\label{Pro1}
\end{align}
where the interpolation points $\{x_k,\,k=0,1,\cdots,\,N\}$ are zeros of {  say, the} \zhu{Legendre} polynomials. Our goal is to identify those points $y_j^{(\bt)}$ where $\prescript{RL}{-1}{\mathcal{D}}_{x}^{\bt} u_N,\,\bt\in(0,1)$, the fractional \zhu{Riemann-Liouville} derivative  of the interpolant, superconverges to $\prescript{RL}{-1}{\mathcal{D}}_{x}^{\bt} u$ in the sense that
\begin{align}
N^\af|\prescript{RL}{-1}{\mathcal{D}}_{x}^{\bt} (u-u_n)(y_j^{(\bt)})|\leq C \max_{x\in[-1,1]} |\prescript{RL}{-1}{\mathcal{D}}_{x}^{\bt} (u-u_n)(x)|, \quad \af>0.\label{def_super}
\end{align}

According to the description \eqref{def_super} of superconvergence points originated from \cite{Zhang+2012+SINUM}, if $\bt=1$ and $\{y_j'\}$ are independent of the particular choice of $u$, we say that $\{y_j'\}$ are the  superconvergence points for the first derivative of the interpolant. To step further, here we want to seek out the points $y_j^{(\bt)}$  independent of the function $u$, at which the corresponding $\bt$-th \zhu{Riemann-Liouville} derivative $\prescript{RL}{-1}{\mathcal{D}}_{x}^{\bt} u_N$ superconverges to $\prescript{RL}{-1}{\mathcal{D}}_{x}^{\bt} u$. In other words, for different order of fractional derivative, different superconvergence points are \zhu{needed}.

\subsection{Analysis} \zhu{Considering} polynomial interpolation for \eqref{Pro1}, the key rule for locating the superconvergece points is to analyze the interpolation error.  \zhu{We are now in the position to show the theoretical analysis}.
\zhu{\begin{proposition}\label{prop_legen}
For the  interpolation of \eqref{Pro1} using collocation points as the zeros of  Legendre polynomials
$L_{N+1}(x)$, Legendre-Lobatto polynomials  $(L_{N+1}-L_{N-1})(x)$, and Legendre-Radau (right and left) polynomials
 $(L_{N+1}\pm L_N)(x)$, the $\bt$-th \zhu{Riemann-Liouville} fractional derivative superconverges at  $\{\xi_k^{(\bt)}\},$ which satisfy
 \begin{align}
\prescript{RL}{-1}{\mathcal{D}}_{x}^{\bt}w_{N+1}(\xi_k^{(\bt)})=0,k=0,1,2,\cdots,N,\label{sptssfracplegen}
\end{align}
where $w_{N+1}(x)$ denotes the  aforementioned four sets of polynomials  respectively.
\end{proposition}}
\begin{proof}
Let $u$ be analytic on $I=[-1,1].$ According to \cite{Zhang+2012+SINUM}, $u$ can be analytically extended to $B_\rho$, which is enclosed by an ellipse $E_\rho$ with $\pm1$ as foci and $\rho>1$ as the sum of its semimajor and semiminor:
\begin{align}
E_\rho: \quad z=\frac12(\rho e^{i\theta}+\rho^{-1} e^{-i\theta}),\quad 0\leq \theta\leq 2\pi.
\end{align}
We consider polynomial $u_N\in P_N$ who interpolates $u$ at $N+1$ points $-1\leq x_0<x_1<\cdots<x_N\leq1.$ The error equation is, according to \cite{Davis+1975+NY}, expressed as
\begin{align}
u(x)-u_N(x)=\frac1{2\pi i}\int_{E_\rho}\frac{w_{N+1}(x)}{w_{N+1}(z)}\frac{u(z)}{z-x}dz,\label{erroreq}
\end{align}
where $w_{N+1}(x)=c\prod_{j=0}^N(x-x_j).$

Noticing Lemma \zhu{\ref{fracLeib}}, taking the fractional derivative of the above equation we have
\zhu{\begin{align}\label{chebanal}
&\prescript{RL}{-1}{\mathcal{D}}_{x}^{\bt}[ u(x)- u_N(x)]
=\frac1{2\pi i}\int_{E_\rho} \left[ \sum_{k=0}^{\infty}{\bt \choose k}{\frac{\prescript{RL}{-1}{\mathcal{D}}_{x}^{\bt-k}w_{N+1}(x)}{(z-x)^{k+1}}}    \right] \frac{u(z)}{w_{N+1}(z)}dz.
\end{align}}

 Let us examine the error equation \eqref{chebanal}{ . According} to the analysis in \cite{Zhang+2012+SINUM}, the leading term of the error { is} determined by $\prescript{RL}{-1}{\mathcal{D}}_{x}^{\bt}w_{N+1}(x).$ At the { $N$} special
points $\{\xi_k^{(\bt)},k=0,\cdots,N\}$, which satisfy $\prescript{RL}{-1}{\mathcal{D}}_{x}^{\bt}w_{N+1}(\xi_k^{(\bt)})=0$, we have the remaining terms, which are usually smaller than the
first term in magnitude at least by a factor $N^\af$.  Here we interpolate the function using \zhu{Legendre polynomials, Legendre-Lobatto polynomials, and Legendre-Radau (right and left) polynomials}, thus the corresponding superconvergence points are \zhu{zeros of}  $\prescript{RL}{-1}{\mathcal{D}}_{x}^{\bt}w_{N+1}(x),$ where $w_{N+1}(x)$ denote the  aforementioned four sets of polynomials respectively.
\end{proof}

\zhu{From Lemma \ref{fracLeibC}, the superconvergence points for the Caputo fractional derivative is easily drawn.}
\begin{remark} Similar to the proof in Propsition  \ref{prop_legen}, for the  interpolation of \eqref{Pro1} using collocation points as the zeros of aforementioned four kinds of \zhu{Legendre} polynomials, the $\bt$-th Caputo fractional derivative superconverges at  $\{z_k^{(\bt)}\},$ which satisfy \begin{align}
\prescript{C}{-1}{\mathcal{D}}_{x}^{\bt}w_{N+1}(z_k^{(\bt)})=0,k=\zhu{0},1,2,\cdots,N,\label{sptssfracpp}
\end{align}
where $w_{N+1}(x)$ denotes the  aforementioned four sets of polynomials  respectively.

\end{remark}
\zhu{As an extension, superconvergence results are derived similarly for Chebyshev polynomials in the following.}
\begin{proposition}\label{prop_cheb}
For the  interpolation of \eqref{Pro1} using collocation points as the zeros of  Chebyshev polynomials of the first kind
{  $T_{N+1}(x)$}, Chebyshev-Lobatto polynomials {  $(T_{N+1}-T_{N-1})(x)$}, and Chebyshev-Radau (right and left) polynomials
{  $(T_{N+1}\pm T_N)(x)$}, the $\bt$-th \zhu{Riemann-Liouville}  fractional derivative superconverges at  $\{y_k^{(\bt)}\},$ which satisfy \begin{align}
\prescript{RL}{-1}{\mathcal{D}}_{x}^{\bt}w_{N+1}(y_k^{(\bt)})=0,k=\zhu{0},1,2,\cdots,N,\label{sptssfracp}
\end{align}
where $w_{N+1}(x)$ denotes the  aforementioned four sets of Chebyshev polynomials  respectively.
\end{proposition}

\zhu{In what follows}, we explain how to \zhu{compute}  $\{\xi_k^{(\bt)}\}$ in \eqref{sptssfracplegen}.

%
%


Denote $L_n(x)$, the Legendre polynomial of degree $n$. In spectral methods, we often use the following four combinations. For convenience, express them as generalized Jacobi functions:
$$
L_n(x) = P^{(0,0)}_n(x), \qquad\text{(Legendre polynomial)}
$$
$$
(L_n-L_{n-2})(x) = P^{(-1,-1)}_n(x), \qquad\text{(Lobatto polynomial)}
$$
$$
(L_n+L_{n-1})(x) = P^{(0,-1)}_n(x), \qquad\text{(Left Radau polynomial)}
$$
$$
(L_n-L_{n-1})(x) = P^{(-1,0)}_n(x). \qquad\text{(Right Radau polynomial)}
$$
Using formula {  \cite{LiHuiyuan}}
\begin{equation}\label{li1}
P^{(\alpha,\beta)}_{n-\beta}(x) =
\frac{\Gamma(n+\alpha-\beta+1)\Gamma(n+1)}{\Gamma(n-\beta+1)\Gamma(n+\alpha+1)} \left( \frac{1+x}{2} \right)^{-\beta}P^{(\alpha,-\beta)}_n(x),
\end{equation}
and
\begin{equation}\label{li2}
\prescript{\zhu{RL}}{-1}D^\mu_x P^{(\alpha,\beta)}_{n-\beta}(x) = \frac{\Gamma(n+\alpha+\mu+1)}{2^\mu\Gamma(n+\alpha+1)} P^{(\alpha+\mu,\beta+\mu)}_{n-\beta-\mu}(x),
\end{equation}
we can derive
\begin{equation}\label{legendl}
\prescript{\zhu{RL}}{-1}D^\mu_x P^{(0,0)}_n(x) = \frac{2^{-\mu}\Gamma(n+1)}{\Gamma(n-\mu+1)} \left( \frac{1+x}{2} \right)^{-\mu}P^{(\mu,-\mu)}_n(x),
\end{equation}
\begin{equation}\label{radaul}
\prescript{\zhu{RL}}{-1}D^\mu_x P^{(\alpha,-1)}_n(x)
= \frac{(n+\alpha)\Gamma(n)}{2^\mu\Gamma(n-\mu+1)} \left( \frac{1+x}{2} \right)^{1-\mu}P^{(\alpha+\mu,1-\mu)}_{n-1}(x),
\quad \alpha=-1, 0.
\end{equation}
Similarly, we have
\begin{equation}\label{legendr}
\prescript{\zhu{RL}}{x}D^\mu_1 P^{(0,0)}_n(x) = \frac{2^{-\mu}\Gamma(n+1)}{\Gamma(n-\mu+1)} \left( \frac{1-x}{2} \right)^{-\mu}P^{(-\mu,\mu)}_n(x),
\end{equation}
\begin{equation}\label{radaur}
\prescript{\zhu{RL}}{x}D^\mu_1 P^{(-1,\alpha)}_n(x)
 = \frac{(n+\alpha)\Gamma(n)}{2^\mu\Gamma(n-\mu+1)} \left( \frac{1-x}{2} \right)^{1-\mu}P^{(1-\mu,\alpha+\mu)}_{n-1}(x),
\quad \alpha=-1, 0.
\end{equation}
Therefore, we conclude that (\ref{sptssfracplegen}) is equivalent to find zeros of:

(1) $P^{(\mu,-\mu)}_{N+1}(x)$ for the Legenrde interpolation, i.e., interpolation at roots of $P^{(0,0)}_{N+1}(x)$;

(2) $(1+x)^{1-\mu}P^{(\mu-1,1-\mu)}_N(x)$ for the Lobatto interpolation, i.e., interpolation at roots of $P^{(-1,-1)}_{N+1}(x)$;

(3) $(1+x)^{1-\mu}P^{(\mu,1-\mu)}_N(x)$ for the Left-Radau interpolation, i.e., interpolation at roots of $P^{(0,-1)}_{N+1}(x)$.

Similar results can be obtained for the right fractional derivative $\prescript{\zhu{RL}}{x}D^\mu_1$.

We see that a fractional derivative of a polynomial of degree $N+1$ may have $N+1$ roots (comparing with $N$ roots for an integer derivative).


\subsection{Numerical validation} \zhu{We report on the numerical results} using the  aforementioned {  four sets of \zhu{Gauss} interpolation points in \eqref{Pro1}. }
%

{\bf Example 3.1 } We choose $u(x)=\frac1{100}(x+1)^{10.15}$ in \eqref{Pro1}.

\begin{figure}[ht!] 
 \centering
\includegraphics[width=4.65in]{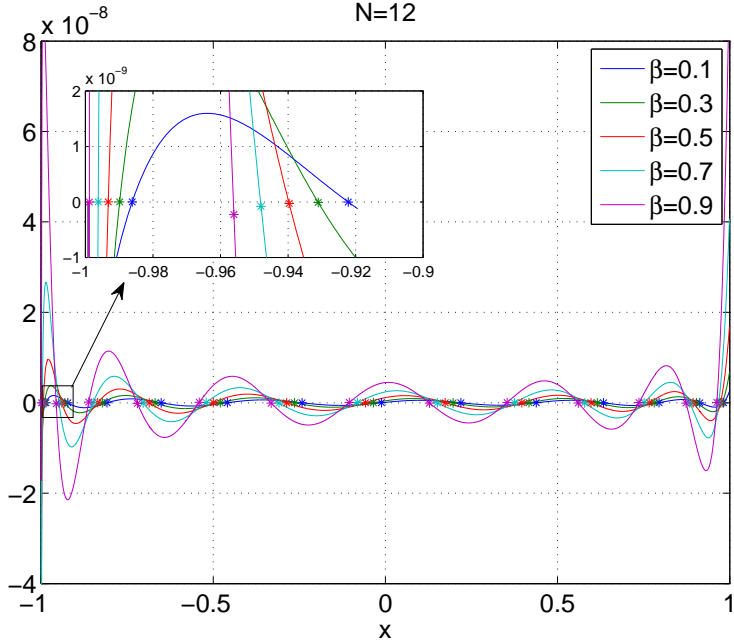}
   \caption{\zhu{Curves of $\prescript{RL}{-1}{\mathcal{D}}_{x}^{\bt}( u- u_N)$} for Example 3.1 using \zhu{Gauss} points, where $*$  denotes the corresponding superconvergence points.}  \label{fig21}
\end{figure}

\begin{figure}[!ht] 
\begin{minipage}{0.99\linewidth}
 \centering
\includegraphics[width=4.65in]{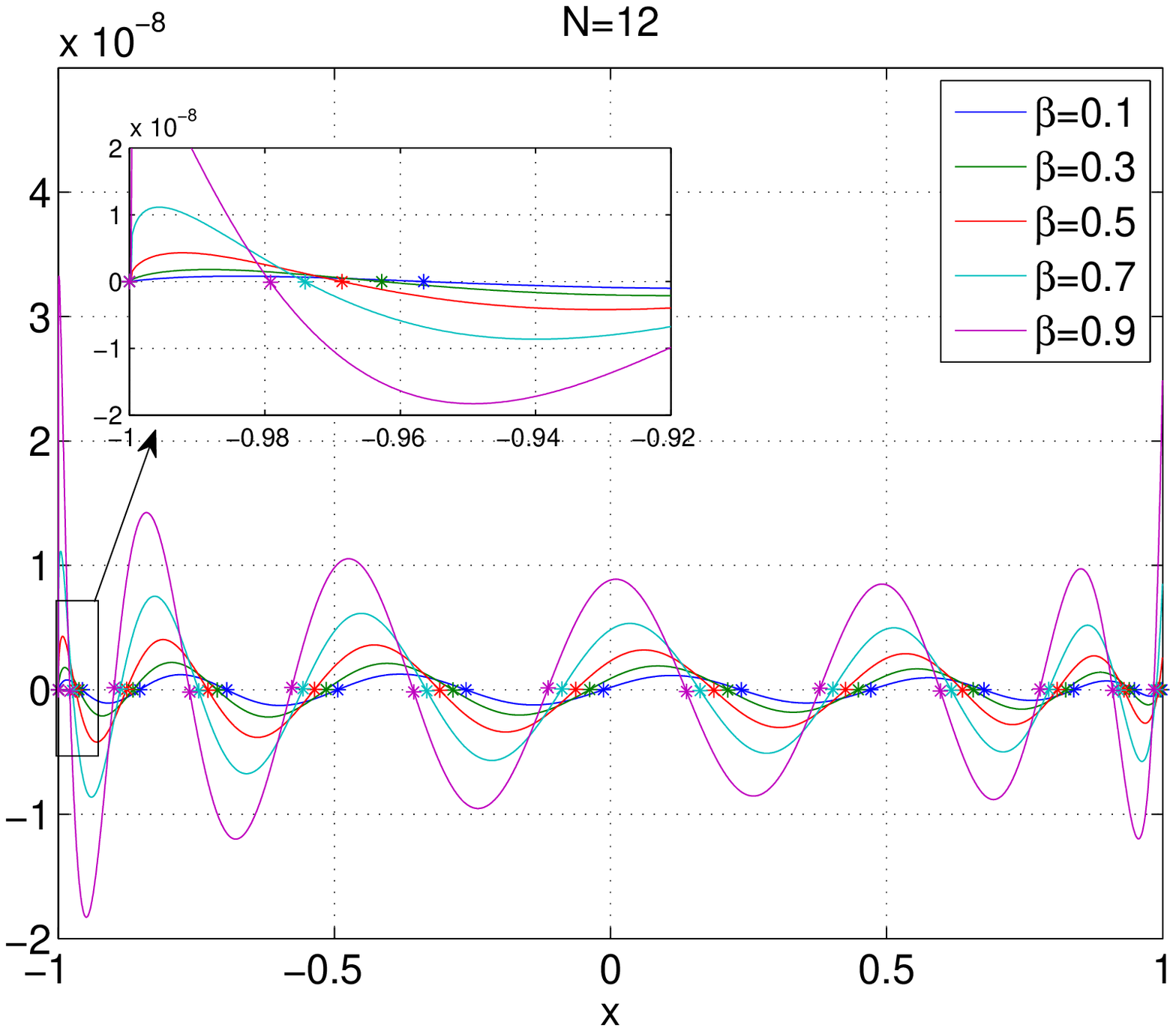}
   \caption{\zhu{Curves of $\prescript{RL}{-1}{\mathcal{D}}_{x}^{\bt}( u- u_N)$}  for Example 3.1 using \zhu{Gauss}-Lobatto points, where $*$  denotes the corresponding superconvergence points.}\label{fig22}
  \end{minipage}
  \begin{minipage}{0.99\linewidth}
 \centering
\includegraphics[width=4.65in]{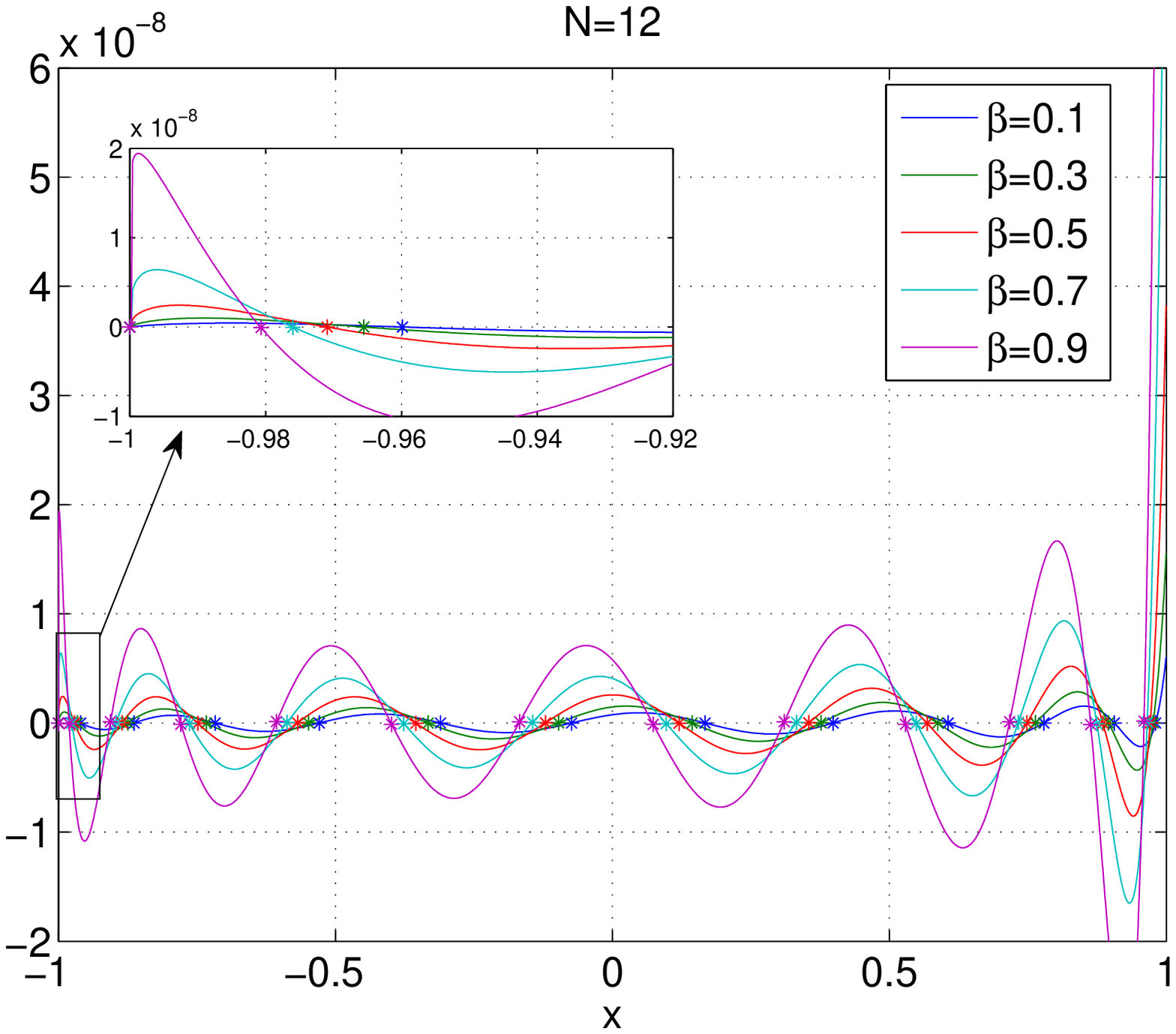}
   \caption{\zhu{Curves of $\prescript{RL}{-1}{\mathcal{D}}_{x}^{\bt}( u- u_N)$}  for Example 3.1 using left \zhu{Gauss}-Radau points, where $*$  denotes the corresponding superconvergence points.}\label{fig23}
  \end{minipage}
\end{figure}

Fractional derivative errors of its interpolant, using \zhu{Gauss} points,   \zhu{Gauss}-Lobatto points, and the left \zhu{Gauss}-Radau points with number of points of $N=12$ are depicted in {  Figures} \ref{fig21}-\ref{fig23} respectively. \zhu{From Proposition \ref{prop_legen}, we know that for a polynomial interpolant of degree {  $N$}, the leading error term is a polynomial of degree { $N+1$}, and there are { $N+1$} roots for its fractional derivative. From {  Figures} \ref{fig21}-\ref{fig23}, it is clearly seen that there are 13 superconvergence points. Particularly, { $-1$} is a superconvergece points for the Gauss-Lobatto and left Gauss-Radau interpolants. To view those roots near { $-1$}, we add zoom-in windows into these figures.}

 \zhu{Five different colors of asterisks denoting the superconvergence points, with orders $\bt=0.1,\,0.3,\,0.5,\,0.7,\,0.9$ for \zhu{Legendre} polynomial interpolation, are show in Figure \ref{fig21}. }It is easily shown that for different orders of fractional derivatives,  the errors at the superconvergence points are significantly smaller than the maximum error.
In addition,  \zhu{it is also indicated that} when the fractional order becomes smaller, the error curve goes lower.

\zhu{Similar situations is easily drawn from  Figure \ref{fig22}-\ref{fig23} for \zhu{Legendre}-Lobatto polynomial and the left \zhu{Legendre}-Radau polynomial respectively.} While for approximations using \zhu{Gauss} points, the errors on the boundaries are larger than that using \zhu{Gauss}-Lobatto points and  left \zhu{Gauss}-Radau points.

\noindent  

\noindent  

\section{Fractional derivative interpolation} In this section, we present the interpolation for the fractional \zhu{Riemann-Liouville} derivatives of a smooth function using GJFs as basis functions and investigate superconvergence points for the function value approximation and fractional derivatives.
{  In particular}, we construct polynomial $u_N$ such that
\begin{align}
\prescript{RL}{x}{\mathcal{D}}_{1}^{s} u_N(x_k)&=\prescript{RL}{x}{\mathcal{D}}_{1}^{s}u(x_k)=f(x_k),\quad u_N(1)=u(1)=0,\nonumber \\
&\quad k=1,2,\cdots,N;\quad-1\leq x_1<\cdots<x_N\leq1,\label{der_interP}
\end{align}
which is equivalent to find the solution for the following factional initial value problem of order $s\in(0,1),$
\begin{align}
\prescript{RL}{x}{\mathcal{D}}_{1}^{s}u(x)=f(x),\quad u{(1)}=0, \quad x\in[-1,1].\label{FJFP}
\end{align}

\subsection{Petrov-Galerkin spectral method for \eqref{FJFP}} Following the same methodology used in \cite{Chen+2014+ar} for solving \eqref{FJFP},
to find $u_N\in\prescript{+}{}{\mathcal{F}}_{N}^{(-s,-s)}$ such that {  $u_N(1)=0$ and}
\begin{align}
(\prescript{RL}{x}{\mathcal{D}}_{1}^{s}u_N, { v})=(f,{ v}), \forall {  v} \in \mp_N . \label{PG_test}
\end{align}
We write the numerical solution as an expansion of  GJF basis
\begin{align}
u_N(x)=\sum_{n=0}^N \tilde{u}_n^{(s)} \prescript{+}{}{J}_{n}^{(-s,-s)}(x)\in \prescript{+}{}{\mathcal{F}}_{N}^{(-s,-s)}(\Lam),\label{uspan}
\end{align}

Taking $v_N=P_k^{(0,0)}$ in Eq. \eqref{PG_test} and using the orthogonality of Legendre polynomials, we derive from \eqref{fracGJFr}
that
\begin{align}
\tilde{u}_n^{(s)} =\frac{n!}{\Gamma(n+s+1)}\tilde{f}_n,\quad 0\leq n \leq N,\label{utilde}
\end{align}
where $\{\tilde{f}_n\}$ is the sequence of coefficients for the Legendre expansion of function $f(x)$. Substituting \eqref{utilde} into \eqref{uspan}, we obtain the
numerical solution:
\begin{align}
u_N(x)=\sum_{n=0}^N \frac{n!}{\Gamma(n+s+1)}\tilde{f}_n \prescript{+}{}{J}_{n}^{(-s,-s)}(x).\label{uspan1}
\end{align}

{  We would like to indicate the relation of the spectral collocation method (4.1) and the Petrov-Galerkin method (4.3). If we select collocation points $x_k$ in (4.1) as the Gauss points, i.e., the roots of the Legendre polynomial of degree $N$, and construct an interpolant
$I_Nf \in {\cal P}_N$ such that $(I_Nu)(1)=0$ and
$(I_Nf)(x_k) = f(x_k)$, for $k = 1, 2, \ldots, N$. Denote $w_k$ as weights of the $N$-point Gauss quadrature.
Next, we multiply both sides of (4.1) by $L_j(x_k)w_k$, sum up over $k$, and obtain
\begin{align}
(\prescript{RL}{x}{\mathcal{D}}_{1}^{s}u_N, L_j)=(I_Nf,L_j), \forall j=1,2,\ldots,N-1, \label{quadrat}
\end{align}
\begin{align}
(\prescript{RL}{x}{\mathcal{D}}_{1}^{s}u_N, L_N)^*=(I_Nf,L_N)^*. \label{quadratn}
\end{align}
Here $*$ indicates that the integration is carry out by the numerical quadrature.
Note that the $N$-point Gauss quadrature is exact for polynomials of degree $2N-1$.
We see that (4.1) is ``almost" equivalent to (4.3) in that only one term (\ref{quadratn}) is done by the $N$-point Gauss quadrature.

Now}, we recall the error estimate of fractional derivative in $L^2$-norm as follows.

\begin{theorem}\label{THPGF}(see\cite{Chen+2014+ar})Let $u$ and $u_N$ be the solution of \eqref{FJFP} and \eqref{PG_test}, respectively.  If $f\in C(\bar{\Lambda})$ and $f^{(l)}\in L^2_{w^{(l-1,l-1)}}(\Lambda)$ for all $1\leq l \leq m,$ then we have that for $1\leq m\leq N+1,$
\begin{equation*}
\|\prescript{RL}{x}{\mathcal{D}}_{1}^{s} (u-u_N)\|+\|u-u_N\|\leq  cN^{-m} \|f^{(m)}\|_{\omega^{(m-1,m-1)}},
\end{equation*}
where c is a positive constant independent of u, N and m.
\end{theorem}

The above error estimate of fractional derivative leads to also the approximation error of function $f$ using  Legendre polynomials. Thus, it naturally yields to the following \zhu{superconvergence} points for fractional derivatives.

\begin{proposition}\label{Frac_SP}
For the Petrov-Galerkin spectral method for \eqref{FJFP}, the \zhu{left} fractional \zhu{Riemann-Liouville} derivative of numerical solution $\prescript{RL}{x}{\mathcal{D}}_{1}^{s} u_N$ superconverges to $\prescript{RL}{x}{\mathcal{D}}_{1}^{s} u$  at the zeros of Legendre polynomial $P_{N+1}^{(0,0)}(x),$ \zhu{namely}, the points $\{\zeta_k^{(0)}\}$ satisfying
\begin{align*}
P_{N+1}^{(0,0)}(\zeta_k^{(0)})=0,\quad k=1,\,2,\cdots,N+1,
\end{align*}
which \zhu{are} usually called Gauss points.
\end{proposition}

{  Due to the special design of basis functions based on different fractional derivatives, the resulting superconvergence points are the same Gauss points for all cases.}

From \eqref{utilde} we see that the resulting system is diagonal, which indicates the direct expansion of the numerical solution using $\left\{\prescript{+}{}{J}_{n}^{(-s,-s)}(x)\right\}$. Therefore, the superconvergence phenomena can also be observed as follows.

\begin{proposition}\label{Frac_SP1}
For the Petrov-Galerkin spectral method for \eqref{FJFP}, the numerical solution $u_N$ superconverges to the exact solution at the zeros of $\prescript{+}{}{J}_{N+1}^{(-s,-s)}(x),$ \zhu{namely}, the points $\{\eta_k^{(s)}\}$  satisfying
\begin{align*}
\prescript{+}{}{J}_{N+1}^{(-s,-s)}(\eta_k^{(s)})=0,\quad k=1,\,2,\cdots,N+1.
\end{align*}
\end{proposition}

In the above proposition, superconvergence points $\{\eta_k^{(s)}\}$,  for the function value, depends on the order of fractional derivative. \zhu{On }the contrary, Proposition \ref{Frac_SP} tells us that the fractional derivative of the approximation superconverges at Gauss points, which are $s$- independent, for any $0<s<1.$ The properties of the two kinds of superconvergence points manifest the essential differences of the estimates in Theorem  \ref{THPGF}.

The following remark gives us the opportunity to see further how the GJFs of \eqref{GJFsr} and \eqref{GJFsl} becoming excellent candidates for solving the fractional equation with underlying singularities at the boundaries.

\begin{remark} For solving problem \eqref{FJFP}, GJF-Petrov-Galerkin spectral method with  \eqref{GJFsr} as basis functions are applied. While, for the following initial value problem \zhu{involving left Riemann-Liouville fractional derivative:}
\begin{align}
\prescript{RL}{-1}{\mathcal{D}}_{x}^{s}u(x)=f(x),\quad u{(-1)}=0, \quad x\in[-1,1].\label{FJFPl}
\end{align}
 one use \eqref{GJFsl}  as basis to construct Petrov-Galerkin spectral method, and the main results \zhu{of} error estimates and superconvergence points could also be proved.
\end{remark}

In the last part of this section, we present numerical examples to confirm the theoretical results and show the superconvergence points for the function value and fractional derivatives.

\subsection{Numerical Examples} In this subsection, we report on numerical results for problem \eqref{FJFP} with three different cases of $f.$
\vskip 0.3cm
{\bf Example 4.1} Take $f(x)=1+x+\cos(x)+\sin(x)$ in \eqref{FJFP}.
\vskip 0.3cm
{\bf Example 4.2} Take $f(x)=e^{\sin(x)+2}$ in \eqref{FJFP}.
\vskip 0.3cm
{\bf Example 4.3} Take $f(x)=(1+x)^{7.89}$ in \eqref{FJFP}.

From Theorem  \ref{THPGF}, we know that the errors of function value and its fractional derivatives decay exponentially regardless of the unknown solution is singular at boundaries or not. Here, we focus on displaying the errors and the superconvergence points for these two cases.

\begin{figure}[ht!] 
 \centering
  \includegraphics[width=4.65in]{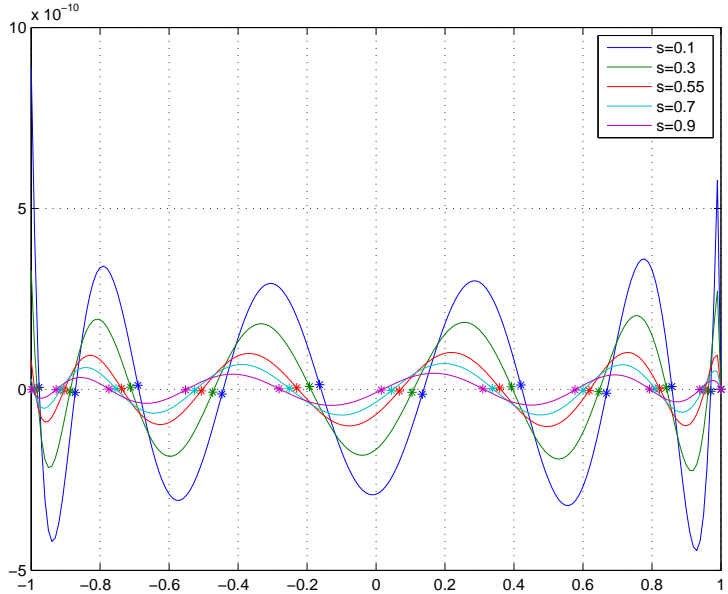}
   \caption{Curves of $u-u_N$ with N=9 for Example 4.1.}\label{fig31}
\end{figure}

\begin{figure}[ht!] 
 \centering
  \includegraphics[width=4.65in]{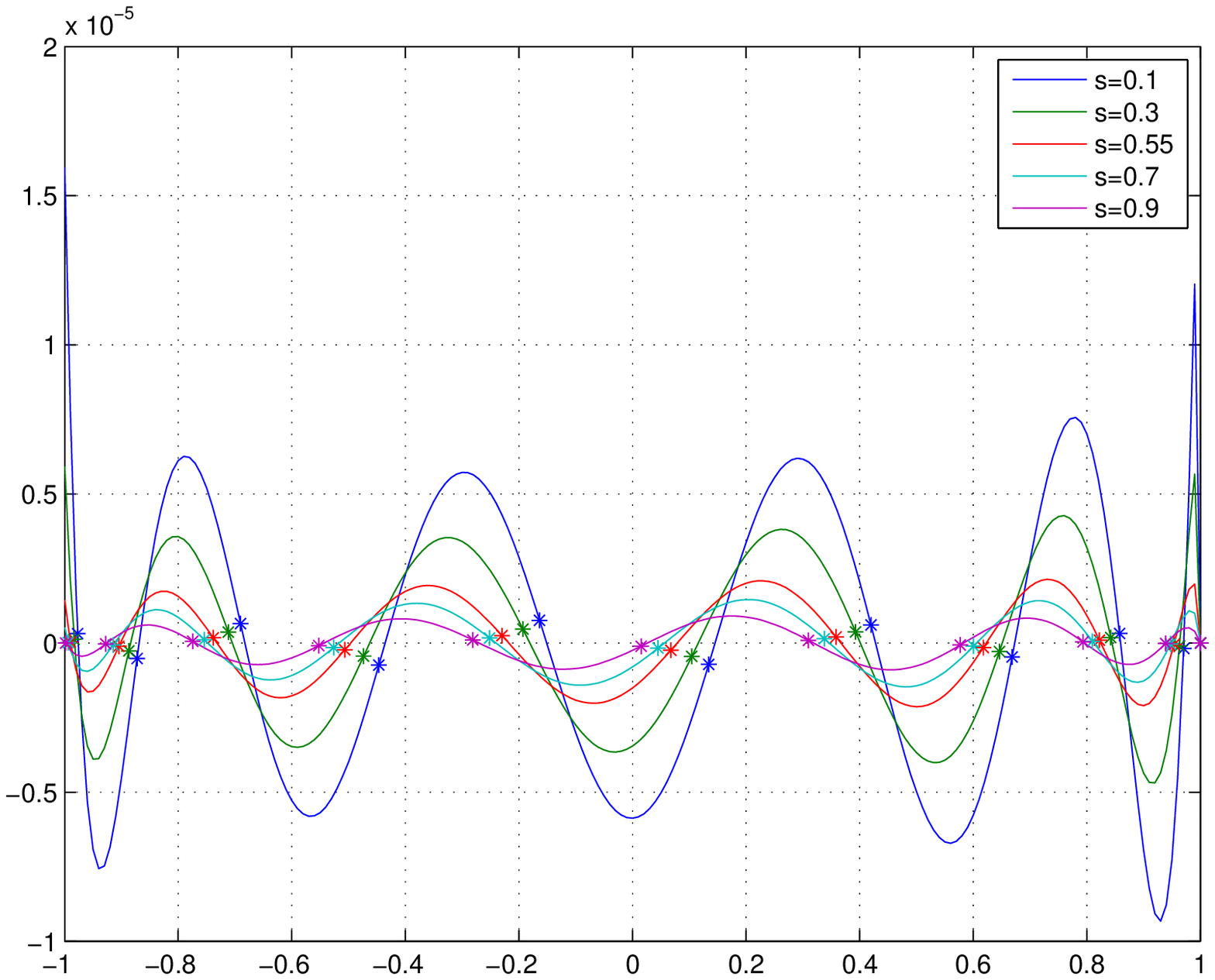}
   \caption{Curves of $u-u_N$ with N=9 for Example 4.2.}\label{fig32}
\end{figure}

\begin{figure}[ht!] 
 \centering
  \includegraphics[width=4.65in]{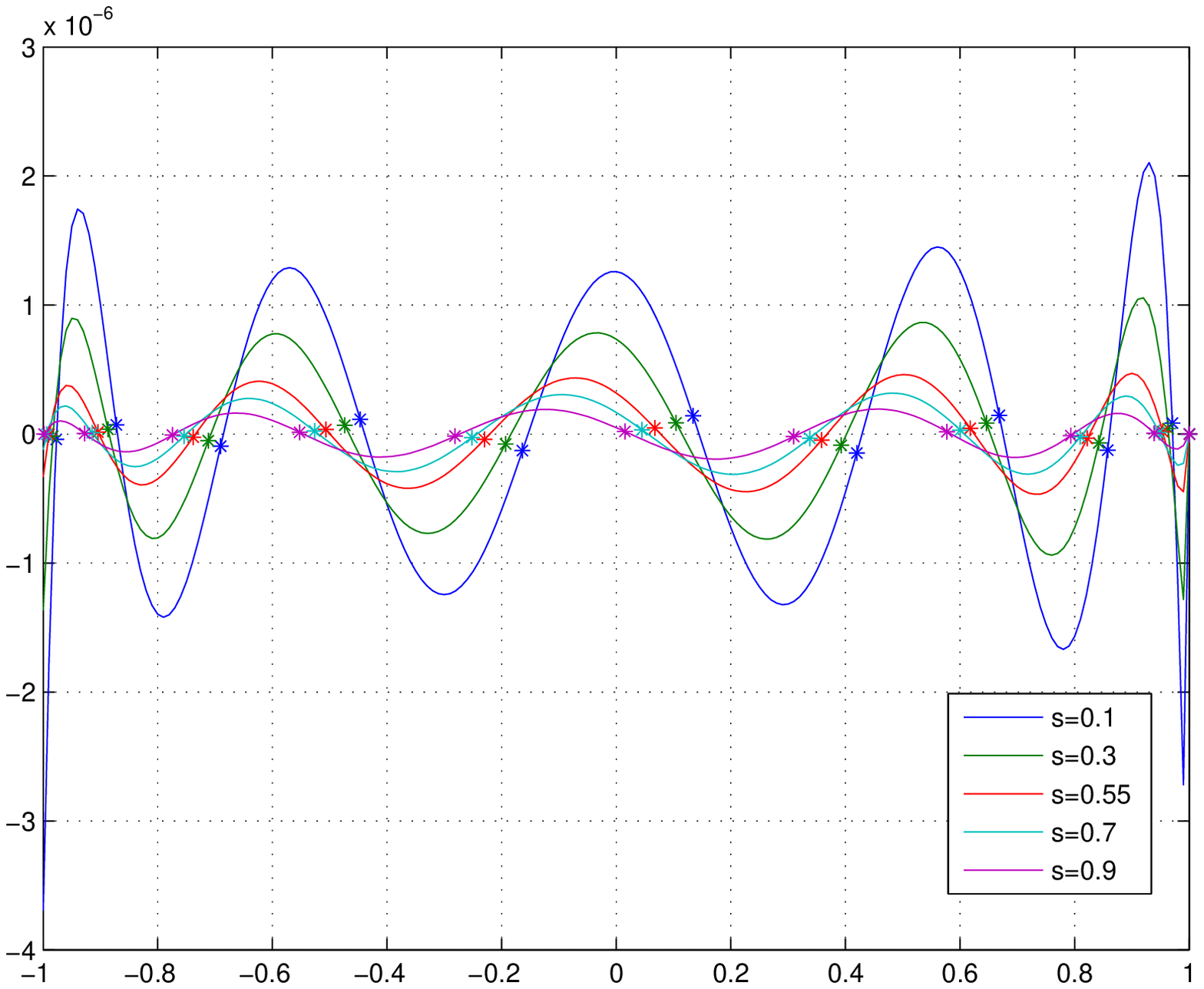}
   \caption{Curves of $u-u_N$ with N=9 for Example 4.3.}\label{fig33}
\end{figure}

We compute a reference exact solution with $N=41$ in \eqref{uspan1}.
  Maximum errors between the exact solution and the numerical solution with $N=9$ for the above examples are {  demonstrated in Figures} \ref{fig31}-\ref{fig33}.  For different orders of fractional derivatives as $s=0.1,0.3,0.55,0.7,0.9,$ the superconvergence points are marked in different colors of asterisks respectively. It is easily seen that  the errors at these points are much smaller than the global maximum errors.
An observation worth noting is that when the order $s\rightarrow1,$ the amplitudes of the error curves getting smaller, each series of the superconvergence points seems converge \zhu{to the standard case of $s=1$}.

\begin{figure}[ht!]
 \centering
  \includegraphics[width=4.65in]{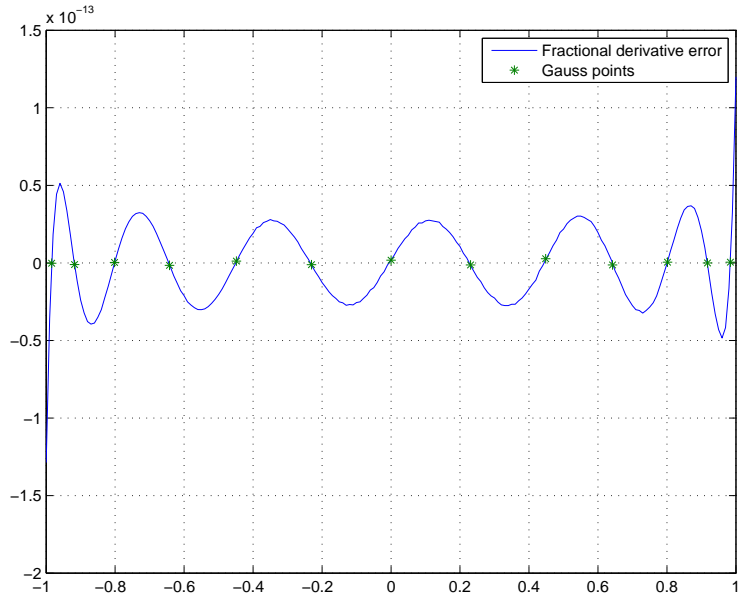}
   \caption{Curves of $\prescript{RL}{x}{\mathcal{D}}_{1}^{s} (u-u_N)$ with N=12 for Example 4.1.}\label{fig31_fra}
\end{figure}

\begin{figure}[ht!]
 \centering
  \includegraphics[width=4.65in]{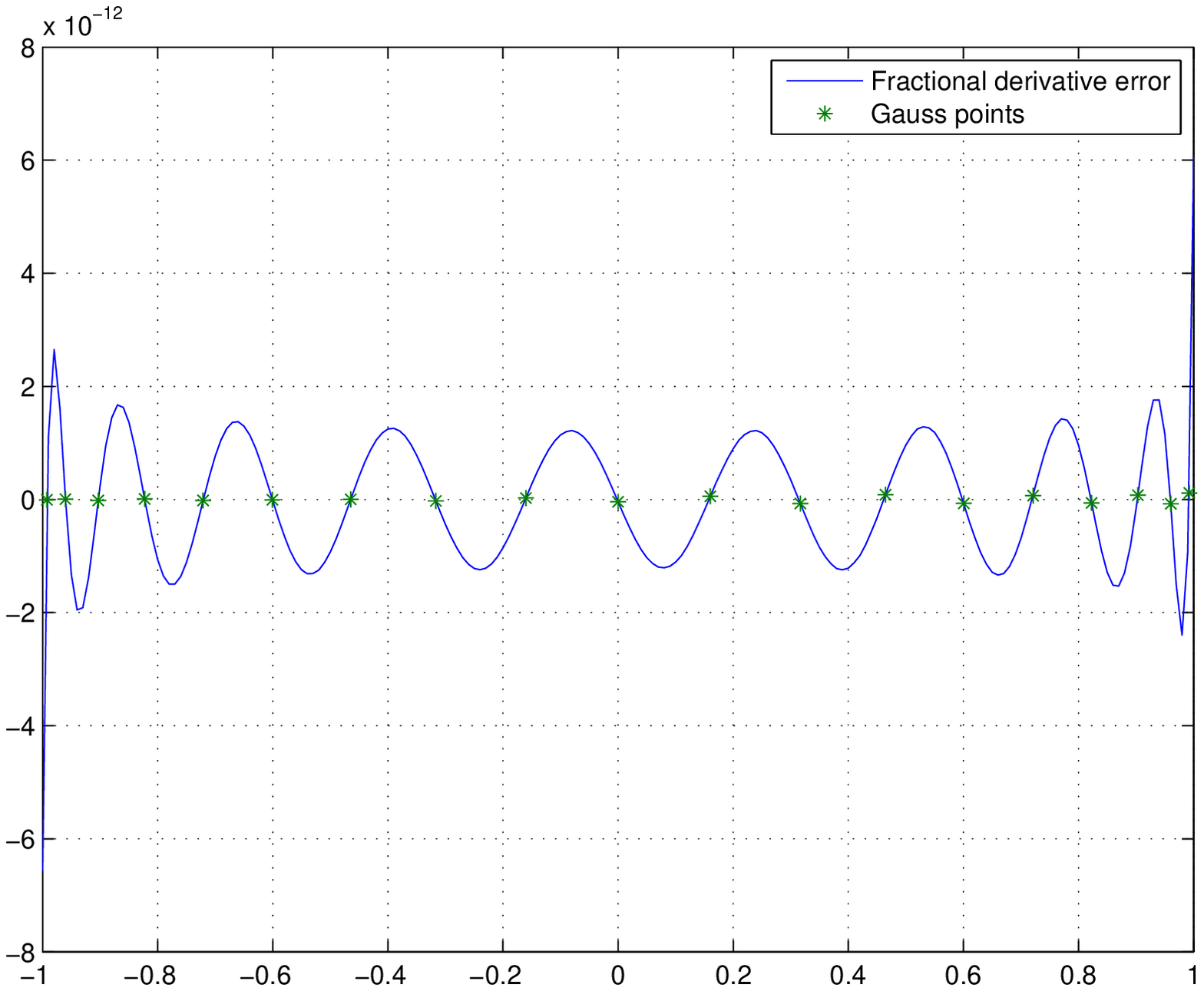}
   \caption{Curves of $\prescript{RL}{x}{\mathcal{D}}_{1}^{s} (u-u_N)$ with N=18 for Example 4.2.}\label{fig32_fra}
\end{figure}

\begin{figure}[ht!]
 \centering
  \includegraphics[width=4.65in]{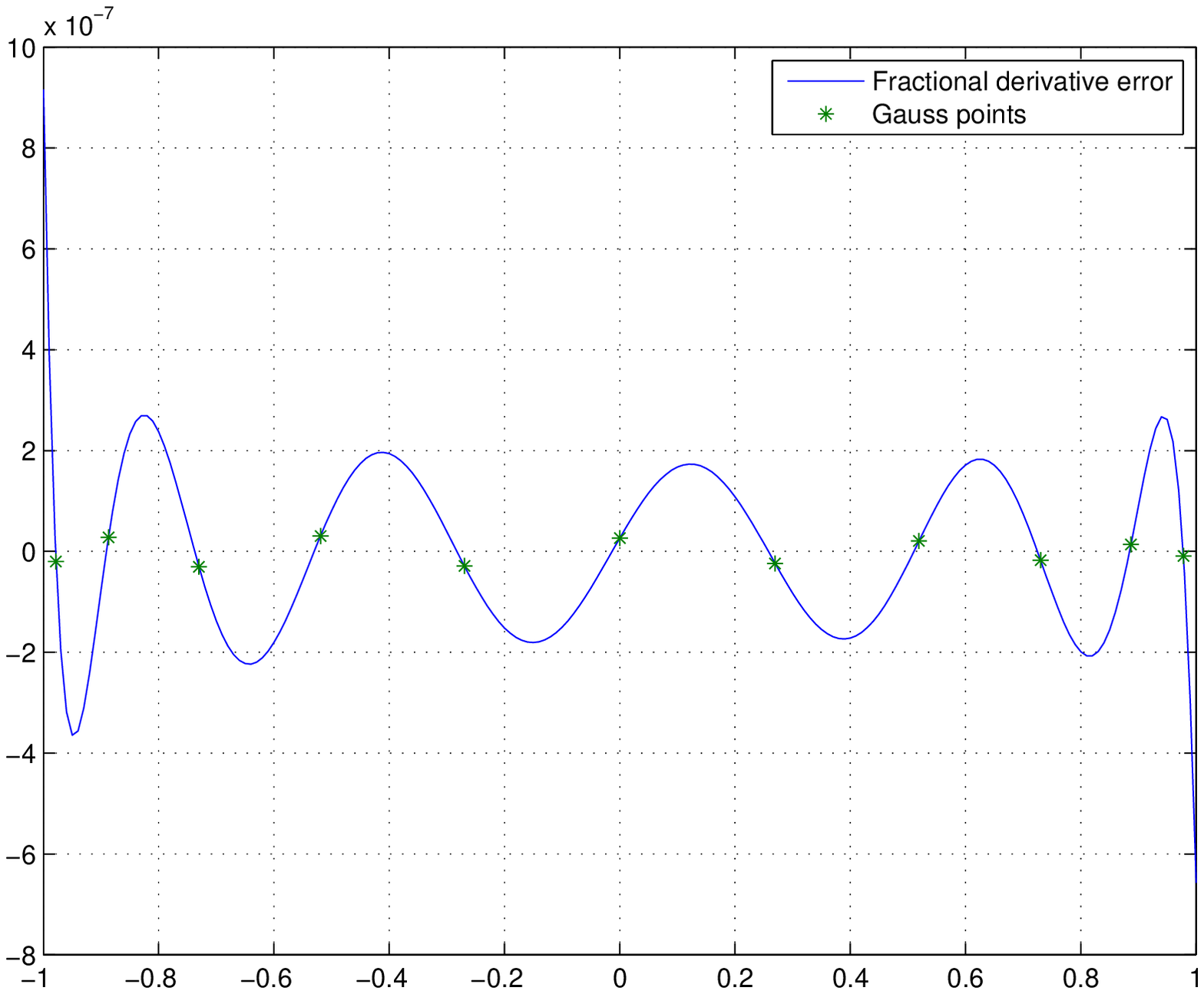}
   \caption{Curves of $\prescript{RL}{x}{\mathcal{D}}_{1}^{s} (u-u_N)$ with N=10 for Example 4.3.}\label{fig33_fra}
\end{figure}

 For the purpose of comparison, we present fractional {  derivative errors in Figures} \ref{fig31_fra}-\ref{fig33_fra} for three examples.
 In contrast {  to the function value approximation, superconvergence} points are the same for different orders of fractional \zhu{derivatives}, which confirms the {  conclusion drawn} in Proposition \ref{Frac_SP}.

\begin{figure}[ht!]
 \centering
  \includegraphics[width=4.65in]{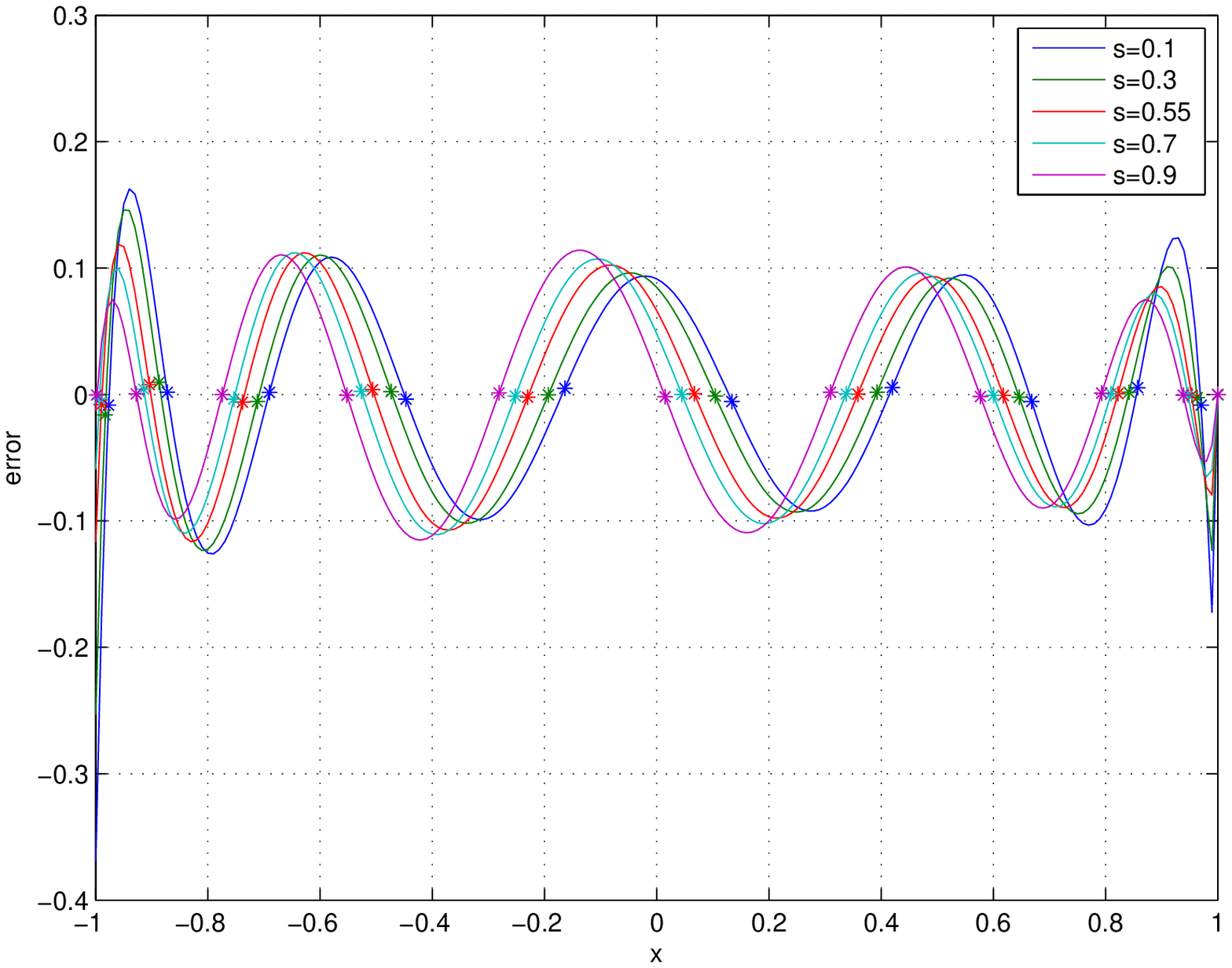}
   \caption{Curves of $u-u_N$ with N=9 for Eq. \eqref{FIVPu}}\label{con}
\end{figure}

{  Remark 4.5. We would like to indicate that although our analysis is for a very special case (4.2), our findings of superconvergence points are actually valid for more general situations as evidence by the following example.
Consider
\begin{align}
\prescript{RL}{x}{\mathcal{D}}_{1}^{s}u(x)+u(x)=f(x),\quad u{(1)}=0, \quad x\in[-1,1]. \label{FIVPu}
\end{align}
We apply the Petrov-Galerkin spectral method (as described above) with $N=9$ for the exact solution $u(x) = (1-x)^{12+\alpha}$. The error curves of function values for different orders of fractional derivative are displayed in Figure \ref{con}, where stars mark the same points as in Proposition \ref{Frac_SP1}. We notice that superconvergence points are the same as for problem (\ref{FIVPu}).}

\section{Concluding remarks}

In this paper, we discuss {  superconvergence phenomena for two kinds of spectral interpolations involving fractional} derivatives. The intended application is the development of high-order {  methods} for fractional problems using spectral methods.

 {  When interpolating function values, superconvergence points for the \zhu{Legendre} basis are located by error analysis.
 This analysis unifies the identification of superconvergence points for fractional and first-order derivatives. When interpolating the fractional derivative values, we found significant differences between superconvergence points for function values and fractional derivatives.
Numerical comparisons between the two cases are provided.

In the future, we plan to study more general fractional differential equations such as (\ref{FIVPu}).}

\vskip.2in

{  {\bf Acknowledgement}. The authors would like to thank Professor Huiyuan Li for introducing formula (\ref{li1}) and (\ref{li2}),
which lead to (\ref{legendl})-(\ref{radaur}).}

\end{document}